\documentclass[12pt]{article}
\usepackage[utf8]{inputenc}
\usepackage{graphicx,psfrag,epsf}
\usepackage{enumerate}
\usepackage{natbib}
\usepackage[margin = 1in]{geometry}
\usepackage{color}
\usepackage{amsmath,amssymb}
\usepackage{bm}
\usepackage{listings}
\usepackage[svgnames]{xcolor}
\usepackage{natbib}

\newcommand{\sumonn}{\sum_{i=1}^n}

\newcommand{\bcon}{\hat{\beta}_{\psi}}
\newcommand{\con}{\Lambda}

\def\T{{ \mathrm{\scriptscriptstyle T} }}

\renewenvironment{abstract}
 {\small
  \begin{center}
  \bfseries \abstractname\vspace{-0mm}\vspace{0pt}
  \end{center}
  \list{}{
    \setlength{\leftmargin}{1cm}%
    \setlength{\rightmargin}{\leftmargin}%
  }%
  \item\relax}
 {\endlist}

\def\spacingset#1{\renewcommand{\baselinestretch}%
{#1}\small\normalsize} \spacingset{1}
\def\T{{ \mathrm{\scriptscriptstyle T} }}

\begin{document}
\title{\bf \LARGE A New Look at $F$-tests}
\date{\vspace{-5ex}}
\maketitle
\vspace{-9mm}
\begin{center}
\textbf{\normalsize Andrew McCormack\textsuperscript{a}, Nancy Reid\textsuperscript{a}, Nicola Sartori\textsuperscript{b}, Sri-Amirthan Theivendran\textsuperscript{a}}
\end{center}
\begin{center}
\textit{\small \textsuperscript{a}Department of Statistical Sciences, University of Toronto\\ \textsuperscript{b}Dipartimento di Scienze Statistiche, Universit{\'a} di Padova}
\\
\vspace{3mm}
\date{\today}
\end{center}
\vspace{3mm}
\begin{abstract}
\noindent Directional inference for vector parameters based on higher order approximations in likelihood inference has recently been developed in the literature. Here we explore  examples of directional inference where the calculations can be simplified, and find that in several classical situations the directional test is equivalent to the usual $F$-test. 
\end{abstract}
\hspace{1cm}\small\textbf{Keywords}: directional tests; exponential family; $F$-test; linear regression

\section{Introduction}
\label{sec:intro}

\subsection{Models and testing}
\label{sec:intro.1}

In many statistical settings we are interested in hypotheses about vector parameters. Examples include testing sets of dummy variables indexing levels of a factor in a regression model, testing interactions in loglinear models for multi-way contingency tables, and tests in models for multivariate responses, such as testing hypotheses for the mean vector or the covariance matrix of a multivariate normal distribution. 
\\
\\
To fix notation, assume we have a model for a response $y_i$ with parametric density function $f_i(y_i;\theta)$. We write $f(y;\theta)$ for the joint density for a sample $y = (y_1, \dots, y_n)$; the maximum likelihood estimator from this sample is $\hat\theta = \hat\theta(y) = \arg\sup_\theta f(y;\theta)$. 
\\
\\
A general approach to testing a given $\theta$ value  is to construct the quadratic form \\ $q(\theta) = (\hat\theta -\theta)^{\T}V^{-1}(\hat\theta - \theta)$, where $V$ is an estimate of the covariance matrix of  $\hat\theta$. Under conditions ensuring that $\hat\theta$ is consistent and asymptotically normally distributed, and that $V$ is a consistent estimator of the covariance matrix of $\hat\theta$, the resulting statistic has a limiting $\chi^2_p$ distribution as $n \rightarrow\infty$, where $p$ is the dimension of $\theta$. An asymptotically equivalent test for $\theta$ is that based on  the log-likelihood ratio statistic
\begin{equation}\label{lrt}
w(\theta) = 2\{\log f(y;\hat\theta) - \log f(y;\theta)\};
\end{equation}	

this also has a limiting $\chi^2_p$ distribution, but the distribution of $q$ and $w$ will be different in finite samples.
Tests about a hypothesis for a subvector of $\theta$ are similarly constructed. Suppose $\theta = (\psi, \lambda)$ where $\psi$, of dimension $d$, is the parameter of interest. The analogous statistics for testing a value $\psi$ are
\begin{eqnarray}
q(\psi) &=& (\hat\psi - \psi)^{\T} V_1^{-1}(\hat\psi - \psi), \label{wald}\\
w(\psi) &=& 2\{\log f(y;\hat\theta) - \log f(y;\hat\theta_\psi)\} \label{lrt2},
\end{eqnarray}
where $V_1$ is an estimate of the covariance matrix of $\hat\psi$, $\hat\theta_\psi = (\psi, \hat\lambda_\psi)$, and $\hat\lambda_\psi$ is the constrained maximum likelihood estimator obtained by maximizing $f(y;\theta)$ over $\lambda$ with $\psi$ fixed; these have under regularity conditions a limiting $\chi^2_d$ distribution. In some contexts a quadratic form based on the score statistic is more convenient for calculations, and has the same limiting distribution.

\subsection{Directional testing}
\label{sec:intro.2}

In the context of linear regression, \citet{Fraser.Massam:1985} proposed  tests that measure departure of $\hat\theta$ from $\theta$, or $\hat\psi$ from $\psi$, in a particular direction on the parameter space. \citet{Skovgaard:1988} considered directional tests in exponential family models, and derived a saddlepoint expansion for the tail probability, based on the $\chi^2$ distribution.   More recently \citet{Davison.etal:2014} and \citet{Fraser.etal:2016} showed how the saddlepoint approximation could be used to compute directional $p$-values via one-dimensional integrals, and illustrated this in a number of models. In this work we present examples for which the directional tests simplify to very well-known omnibus tests. This sheds light on directional testing and gives a new look at some  familiar test statistics. 
\\
\\
The theory and methods for directional tests are given in \citet{Davison.etal:2014} and \citet{Fraser.etal:2016}, and provided for completeness in the Supporting Information. Suppose we have an exponential family model with sufficient statistic $u=u(y)$ and canonical parameter $\varphi$:
\begin{equation}\label{expfam}
f(y;\varphi) = \exp\{\varphi^{\T}u - \kappa(\varphi)\}h(y),
\end{equation}
and we are interested in the hypothesis $H_\psi:\psi(\varphi) = \psi$.
We write $\hat\varphi_\psi$ for the constrained maximum likelihood estimate of $\varphi$ under $H_\psi$. The directional test of $H_\psi$ restricts attention to the line in the sample space joining $u_{\psi}$ with $u^0$, where $u_{\psi}$ is the value of the sufficient statistic which would give the $\hat\varphi_\psi$ as the maximum likelihood estimate, and $u^0$ is the observed value of the sufficient statistic.  The directional $p$-value is the tail area probability for the length $||u_\psi - u^0||$, conditional on the direction $(u_\psi-u^0)/||u_\psi - u^0||$. It is sometimes convenient to define $s = u - u^0$, so that $s^0 = 0$ and $s_\psi = u_\psi-u^0$.
\\
\\
The $p$-value for this directional approach is defined by a ratio of two integrals:
\begin{equation}\label{pval}
p(\psi) = \frac{\int_1^{t_{max}}t^{d-1}h(t;\psi) \, dt}{\int_0^{t_{max}}t^{d-1}h(t;\psi) \, dt},
\end{equation}
where $d$ is the dimension of $\psi$, and $t$ indexes points along the line, with $t=0$ corresponding to the value $u_\psi$, and $t=1$ corresponding to the observed value $u^0$. The two one-dimensional integrals in (\ref{pval}) can be easily and accurately computed numerically. The factor $t^{d-1}$ comes from the Jacobian when changing to spherical coordinates on the plane.
The ingredients needed for the calculation of (\ref{pval}) are the log-likelihood function,  maximum likelihood estimate, and the observed Fisher information as functions of $t$, as well as the observed value of the constrained maximum likelihood estimate. An expression for $h(t;\psi)$ when the hypothesis is linear in $\varphi$ is given in (8) in \citet{Davison.etal:2014} and when the hypothesis is nonlinear in \citet{Fraser.etal:2016}. For completeness the general expression for $h(t;\psi)$ is presented here and described in detail in the Supporting Information:
\begin{equation}\label{densline}
h(t;\psi) \propto \exp\left[\ell\{\hat\varphi_\psi;s(t)\}-\ell\{\hat\varphi;s(t)\}\right]|\hat J_{\varphi\varphi}|^{-1/2} |\tilde J_{(\lambda\lambda)}|^{1/2},
\end{equation}
where $\ell(\varphi;s) = \varphi^{\T}(\theta)s+\ell^0\{\theta(\varphi)\}$,
$\ell^0(\theta)$ is the observed value of the log-likelihood function from (\ref{expfam}), and the score variable $s$ is constrained to the line $s(t), t > 0$.   The factor  $\vert \tilde J_{(\lambda\lambda)}\vert$, is needed only when the hypothesis is nonlinear in the canonical parameter. The calculation of the information determinants in \eqref{densline} are described in detail in the Supporting Information.
\\
\\
If the underlying model is not an exponential family model, an initial approximation to that model, called the tangent exponential model, is constructed first, and the arguments apply again within this model \citep[Ex.4.3]{Fraser.etal:2016}.
\\
\\
We show in this paper that the ratio of integrals in (\ref{pval}) can be calculated explicitly in some simple models, which helps to explain the accuracy of the approximation as evidenced in \\ \citet{Davison.etal:2014} and \citet{Fraser.etal:2016}. We present the calculation for simple examples first, and then use similar techniques to provide a new view of classical tests in multivariate normal models.  

\section{Inference in scale models}

\subsection{Scalar parameter of interest}

In a one-dimensional sub-model, the directional test gives two-sided $p$-values as it reduces to the probability of the right (or left) tail, conditional on being in that tail.
In the two examples in this section we demonstrate this, as the calculations can be carried out analytically, and the arguments motivate exact calculations in the regression setting of \S 3.1. The results help to explain the accuracy of the simulations in \citet{Davison.etal:2014}. 
 
\subsection{Comparison of Exponential Rates}\label{sec:scale.1}
We consider first the example of \citet[\S 5.2]{Davison.etal:2014} in the case of just two groups ($g=2$). 
Suppose that $y_{ij}$ are independent random variables following an exponential distribution with rates $\theta_{i}$ for $i=1,2$ and $j=1,\dotsc,n_i$, and we would like to test the null hypothesis $H_{\psi}: \theta_1/\theta_2 = \psi$ for some $\psi \in (0,\infty)$. An exact test is available, since $W_\psi = \psi\bar{y}_1/\bar{y}_2$ follows an $F(2n_1, 2n_2)$ distribution under $H_{\psi}$. 
\\
\\
The log-likelihood of the full model is
\begin{align}
\ell(\theta ; y) = \sum_{i = 1}^2\sum_{j = 1}^{n_i} (\log\theta_i - \theta_iy_{ij} ).
\end{align}
The canonical parameter is $\varphi(\theta) = (-\theta_1, -\theta_2)$, the sufficient statistic is $u = (u_1, u_2) = (\Sigma_j y_{1j}, \Sigma_j y_{2j})$,  the maximum likelihood estimate of $\theta$ is $\hat\theta = (n_1/u_1, n_2/u_2)$, and the constrained maximum likelihood estimator is 
\[\hat\theta_\psi = \left(\frac{n\,\psi}{\psi u_1+u_2}, \frac{n}{\psi u_1,+u_2}\right), 
\]
giving
\[
u_\psi = (u_{1\psi}, u_{2\psi})  = \frac{1}{n}\left(u_1n_1 + \frac{u_2n_1}{\psi}, \psi u_1n_2 + u_2n_2\right).
\]

By a standard property of exponential families $u_{\psi}$ is the expected value of $u$ under $H_{\psi}$. We define the line $s(t)$, $t \geq 0$ to be $(1-t)s_{\psi}$ where $s_{\psi} = u_{\psi} - u^0$. At $t = 1$, $s(t)$ is equal to the observed value of $s$, $s^0 = 0$. The saddlepoint approximation to the density of $s$ on the ray $s(t)$ is 
\begin{align}
t^{d-1}h(t;\psi) = (1-t/a_1)^{n_1-1}(1-t/a_2)^{n_2-1},
\label{Expdensity}
\end{align}
with $a_i =  u_{i\psi}/(u_{i\psi} - u_i)$ for $i = 1,2$. The support of the density is over $s(t)> -u^0$, or equivalently $\hat{\theta}\{s(t)\} \geq 0$, so $t_{max}= \max(a_1, a_2)$. 
\\
\\
Although the hypothesis is not linear in the canonical parameter, it turns out that  the nuisance parameter adjustment term $\vert \tilde J_{(\lambda\lambda)}\vert$ is independent of $t$ because $\varphi(\theta)$ is linear in the nuisance parameter $\theta_2$. 
\\
\\
Suppose that $t_{max} = a_2$ so that $\psi \bar y_1 \geq \bar y_2$, then the directional $p$-value is 
\begin{align}
p(\psi) = \frac{\int_1^{a_2} (1-t/a_1)^{n_1-1}(1-t/a_2)^{n_2-1}\, dt}
{\int_0^{a_2} (1-t/a_1)^{n_1-1}(1-t/a_2)^{n_2-1}\, dt}.
\label{pvaluecondit}
\end{align}

Let $p_{\text{num}}$ and $p_{\text{denom}}$  be the numerator and denominator of \eqref{pvaluecondit} respectively. Taking \\ $x = (1-t/a_1)/(1-t/a_2)$ gives
\begin{align*}
p_{\text{num}}
= c\int_{\psi\bar y_1 / \bar y_2 }^\infty
x^{n_1-1}
\left(
1+\frac{n_1}{n_2}x
\right)^{-(n_1+n_2)}
\, dx
,
\quad
p_{\text{denom}}
= c\int_{1}^\infty
x^{n_1-1}
\left(
1+\frac{n_1}{n_2}x
\right)^{-(n_1+n_2)}
\, dx.
\end{align*}
The constant $c$ is the same in $p_{num}$ and $p_{den}$ because the same change of variables is used in both integrals.
Recognizing the above integrands as the density of a $F(2n_1,2n_2)$ random variable with CDF $G(x)$, $\eqref{pvaluecondit}$ becomes 
\[
p(\psi) = \frac{1-G(W_\psi)}{1-G(1)},
\]
showing how the directional $p$-value is related to that based on the $F$-test of $W_\psi =\psi\bar{y}_1/\bar{y}_2$. Similarly, when  $\psi \bar y_1 < \bar y_2$ the directional $p$-value is $p(\psi) = G(W_\psi)/G(1)$, so that  
\begin{align}
p(\psi) = \text{I}(W_\psi \geq 1)\frac{1-G(W_\psi)}{1-G(1)} + \text{I}(W_\psi < 1)\frac{G(W_\psi)}{G(1)}.
\label{FormOfPVal}
\end{align}
The terms that are multiplied by the indicator functions in $\eqref{FormOfPVal}$ are uniformly distributed over $[0,1]$ when conditioned on $W_\psi$ in the appropriate region. It follows that the directional $p$-value is also uniformly distributed since for $x \in [0,1]$
\begin{align*}
P\{p(\psi) \leq x\} & = P\bigg\{\frac{1 - G(W_\psi)}{1-G(1)} \leq x \mid W_\psi \geq 1\bigg\}P(W_\psi \geq 1) + P\bigg\{\frac{G(W_\psi)}{G(1)}\leq x \mid W_\psi < 1\bigg\}P(W_\psi < 1)
\\
& = xP(W_\psi\geq 1) + xP(W_\psi < 1) = x .
\end{align*}
A standard one-tailed $F$-test in this scenario would output $p$-values of $G(W_{\psi})$ if $W_{\psi} < 1$ and $1 - G(W_{\psi})$ otherwise. Such a test does not provide uniformly distributed $p$-values as \\ $p(\psi) \leq \max\{G(1), 1- G(1)\} < 1$. The directional test corrects the $p$-values from the one-tailed $F$-test by the appropriate tail probabilities $G(1)$ or $1 - G(1)$.

  \subsection{Comparison of Normal Variances}

Suppose that $y_{ij}\sim N(\mu_i, \sigma_i^2)$ are independent random variables for $i=1,2$ and $j=1\dotsc,n_i$, and we wish to test $H_{\psi}: \sigma_1^2/\sigma_2^2 = \psi$. Under the  hypothesis, $W_\psi= \psi s_2^2/s_1^2$ follows an $F(n_2-1, n_1-1)$ distribution, where $s_i^2 = n_i(n_i-1)^{-1}v_i^2$ is the unbiased sample variance estimate for group $i$, with $v_i^2 = n_{i}^{-1}\Sigma_{j=1}^{n_i} (y_{ij}-\bar{y}_i)^2$.
\\
\\
Following a derivation similar to that in \citet[\S 5.1]{Davison.etal:2014}, 
the directional $p$-value is
\begin{align}
p(\psi) = \frac{\int_1^{1/a_1} (1-ta_1)^{(n_1-3)/2}(1-ta_2)^{(n_2-3)/2}\, dt}
{\int_0^{1/a_1} (1-ta_1)^{(n_1-3)/2}(1-ta_2)^{(n_2-3)/2}\, dt} ,
\label{pnormcondit}
\end{align}

where  $a_i=(\hat{\sigma}_{i\psi}^2 -v_i^2)/\hat{\sigma}_{i\psi}^2$, $i=1,2$, and $\hat{\sigma}_{i\psi}^2$ is the constrained maximum likelihood estimator for ${\sigma}_i^2$.  
The same change of variable as in (\ref{pvaluecondit}) gives
\[
p(\psi) = \frac{1 - G(W_\psi)}{1 - G\{\frac{n_2(n_1 - 1)}{n_1 (n_2 - 1)}\} },
\]

when $v_{1}^2\le \psi v_{2}^2$, or equivalently, $(n_1 - 1)n_2 / (n_2 - 1)n_1 \le \psi s_2^2 / s_1^2$, where $G(x)$ is the CDF of a $F(n_2-1,n_1-1)$ random variable. Combining this with the case $v_1^2 > \psi v_2^2$ gives 
\begin{align}
p(\psi) 
= I\bigg\{W_\psi \ge \frac{n_2(n_1-1)}{n_1(n_2 - 1)} \bigg\} \frac{1 - G(W_\psi)}{1 - G\{\frac{n_2(n_1-1)}{n_1(n_2 - 1)}\}} + I\bigg\{W_\psi < \frac{n_2(n_1-1)}{n_1(n_2 - 1)} \bigg\}\frac{G(W_\psi)}{G\{\frac{n_2(n_1-1)}{n_1(n_2 - 1)}\}} .
\label{NormVarDirp}
\end{align}

As noted in \citet[Ex. 5.1]{Davison.etal:2014}, this expression simplifies to the two-sided $F$-test if $n_1=n_2$ since $v_1^2 / v_2^2 = s_1^2 / s_2^2$ in this case. When $n_1 \ne n_2$ \eqref{NormVarDirp} is the exact $F$-test based on the tail probabilities of $\psi$ times the ratio of the biased maximum likelihood estimators $v_i^2$.               
With more than two groups there is no exact test for comparison, but simulations in \citet{Davison.etal:2014} show that the directional test is very accurate even with a very large number of groups, and hence large numbers of nuisance parameters, whereas the usual likelihood ratio test (\ref{lrt2}) breaks down.

\label{sec:scale.2}

\section{Inference for location parameters}
\label{sec:location}
\subsection{Linear Regression}

We now consider testing the the null hypothesis $H_{\psi}: A\beta = \psi$ in a linear regression model: 
\[
y_i = x_i^{\T}\beta + \epsilon_i, i = 1,...,n,
\]
where both $X_i$ and $\beta$ are vectors of length $p$, and $\epsilon_i$ are independently distributed as $N(0,\sigma^2)$ with an unknown variance. In the linear constraint $A$ is a given $d\times p$ matrix and  it is assumed that $A$ has maximal rank, so the dimension of the parameter of interest is $d$ and that of the implicit nuisance parameter is $p+1-d$. This null hypothesis encompasses many other hypotheses of interest, such as testing for the equality of group means when the group variances are equal.  If $X$ is taken to be the matrix with rows $x_i^{\T}$, the log-likelihood function for  $\theta^{\T} = (\beta^{\T},\sigma^2)$ is
\begin{align*}
\ell(\theta; y) = -\frac{n}{2}\log\sigma^2 -\frac{1}{2\sigma^2}\big(y^{\T}y - 2y^{\T}X\beta +\beta^{\T}X^{\T}X\beta \big).
\end{align*}
This can be re-expressed in exponential family form with canonical parameter $\varphi(\theta)^{\T} = \sigma^{-2}(\beta^{\T},-1/2)$ and sufficient statistic $u^{\T} = (y^{\T}X,y^{\T}y)= (u_1^{\T}, u_2)$. The unconstrained and constrained maximum likelihood estimates for $\beta$ respectively are
\begin{align*}
\widehat{\beta} = (X^{\T}X)^{-1}X^{\T}y, \quad
\widehat{\beta}_{\psi} =  \widehat{\beta} - (X^{\T}X)^{-1}A^{\T}\big\{ A(X^{\T}X)^{-1}A^{\T}\big\}^{-1}(A\widehat{\beta} - \psi).
\end{align*}
Accordingly, the unconstrained and constrained estimates of $\sigma^2$ are $\hat{\sigma}^2 = n^{-1}\Sigma(y_i-x_i^{\T}\widehat\beta)^2$ and \\ $\hat{\sigma}_{\psi}^2=n^{-1}\Sigma(y_i-x_i^{\T}\widehat\beta_{\psi})^2$. The value of $s$ that has $\hat{\theta}_{\psi}$ as the maximum likelihood estimate of $\theta$ is \begin{align*}
s_{\psi}^{\T} = (\widehat{\beta}_{\psi}^{\T}X^{\T}X - y^{\T}X,n\hat{\sigma}^2_{\psi} + \widehat{\beta}_{\psi}^{\T}X^{\T}X\widehat{\beta}_{\psi} - y^{\T}y);
\end{align*}
On the line $s(t)$, $t \geq 0$, the log-likelihood function for $\varphi(\theta)$ is
\begin{align}
\ell\{\varphi(\theta);s(t)\} = -\frac{n}{2}\log\sigma^2 -\frac{1}{2\sigma^2}\big\{ u_2(t)  - 2u_1(t)\beta + \beta^{\T}X^{\T}X\beta\big\} ,
\label{likehdofstlinconst}
\end{align}
with $\{ u_1(t), u_2(t)\} = u(t) = u^0 + s(t)$. From $\eqref{likehdofstlinconst}$ we obtain the maximum likelihood estimates of $\beta$ and $\sigma^2$ as functions of $t$: $\widehat{\beta}(t) = (X^{\T}X)^{-1}u_1(t)$ and \\ $\hat{\sigma}^2(t) = n^{-1} \big\{u_2(t) - 2u_1(t)\widehat{\beta}(t) + \widehat{\beta}^{\T}(t)X^{\T}X\widehat{\beta}(t) \big\}$. 
\\
\\
Using (\ref{densline}) we get
\begin{align}
h(t;\psi) 
&= \{\hat{\sigma}^2(t)\}^{(n - p - 2)/2} = \big\{ \hat{\sigma}^2_{\psi} - \frac{t^2}{n}(y - X\widehat{\beta}_{\psi})^{\T}X(X^{\T}X)^{-1}X^{\T}(y - X^{\T}\widehat{\beta}_{\psi})  \big\}^{(n - p - 2)/2}.
\label{NormRegresDensity}
\end{align}

The density along the line $s(t)$, passing through $s_{\psi}$ and $s^0$ is 
\begin{align}
t^{d-1}h(t;\psi) = t^{d-1}\big\{ \hat{\sigma}^2_{\psi} - \frac{t^2}{n}(y - X\widehat{\beta}_{\psi})^{\T}X(X^{\T}X)^{-1}X^{\T}(y - X^{\T}\widehat{\beta}_{\psi})  \big\}^{(n - p - 2)/2} = t^{d - 1}(a - bt^2)^{(n - p - 2)/2}.
\end{align}
As the hypothesis can be expressed as a linear function of the canonical parameter, there is no need for the nuisance parameter adjustment term $\vert \tilde J_{(\lambda\lambda)} \vert $. 
\\
\\
To compute the directional $p$-value (\ref{pval}) we need $t_{max}$, which here is the largest value of $t$ for which $ \hat{\sigma}^2(t) \geq 0$:
\begin{align*}
t_{max} = \left\{\frac{n\hat{\sigma}^2_{\psi}}{(y - X\hat{\beta}_{\psi})^{\T}X(X^{\T}X)^{-1}X^{\T}(y - X\hat{\beta}_{\psi})}\right\}^{1/2} = \big(\frac{a}{b}\big)^{-1/2} .
\end{align*}

If we make the change of variables 
\begin{align*}
x = \frac{n - p}{d}\big(\frac{a}{bt^2} - 1\big),
\end{align*}
a computation detailed in the Supporting Information verifies that \eqref{pval} simplifies to
\begin{align}
p(\psi) = 1 - \text{G}\bigg[ \frac{(A\hat{\beta} - \psi)^{\T}\{A(X^{\T}X)^{-1}A^{\T}\}^{-1}(A\hat{\beta} - \psi) }{d} \bigg/ \frac{(y-X\hat{\beta})^{\T}(y-X\hat{\beta})}{n - p}\bigg] ,
\label{linconstpval}
\end{align}

where $G(x)$ is the cumulative distribution function of a $F_{d,n-p}$ random variable. The test statistic in the cumulative distribution function of \eqref{linconstpval} is the standard test statistic for $H_\psi$ and is exactly distributed as a $F_{d,n-p}$ random variable; see for example \citet[Ch. 8]{Rencher and Schaalje:08}.
\\
\\
This result gives a new interpretation of the $F$-statistic: it measures the magnitude of the sufficient statistic for $\psi$, conditional on the direction indicated by the observed data. As the normal distribution is spherically symmetric, the magnitude is distributed independently of the direction.

\subsection{Hotelling's $T^2$}

The result in the previous section suggests comparing the directional test for a multivariate normal mean to Hotelling's $T^2$ statistic.  Suppose $y_i$, $i = 1,...,n$ are independent observations from the multivariate normal distribution, $N_d(\mu,\Lambda^{-1})$, with unknown covariance matrix $\Lambda^{-1}$. The full parameter is $\theta = \{\mu, \text{vec}(\Lambda)\}$, where $\text{vec}$ gives a vectorization of the columns of a matrix.  We seek to test the hypothesis $H_{\psi}: \mu = \psi$. 
The distribution for $y = (y_1, \dots, y_n)$ is an exponential family model, with  canonical parameter $\varphi^{\T}(\theta) =  \{\mu^\T \Lambda , \text{vec}^\T (\Lambda)\}$,  and sufficient statistic $u^{\T} = \{n\bar y^{\T}, \text{vec}^\T(\Sigma y_iy_i^{\T})\}$. The unconstrained maximum likelihood estimates are $\hat{\mu} = \bar{y}$ and $\hat{\Lambda} = n\{\Sigma_{i}(y_i - \bar{y})(y_i - \bar{y})^{\T}\}^{-1}$ while the constrained maximum likelihood estimate for $\Lambda$ is $n\{\Sigma_{i}(y_i - \psi)(y_i - \psi)^{\T}\}^{-1}$. Under $H_{\psi}$, the expected value of the centered sufficient statistic $s$ is
\begin{align}
s^{\T}_{\psi} = \bigg[n\psi - n\bar{y}, \text{vec}^\T\big\{\frac{n}{2}(\psi\bar{y}^{\T} + \bar{y}\psi^{\T} - 2\psi\psi^{\T})\big\}\bigg].
\end{align}
As shown in the Supporting Information the maximum likelihood estimators 
on the ray $s(t)$ are
\begin{align}
\hat{\varphi}^{\T}(t) & = \big[\{\psi + t(\bar{y} - \psi)\}^\T\hat{\Lambda}(t), \text{vec}^\T\{\hat{\Lambda}(t)\}\big],
\\
\hat{\Lambda}^{-1}(t) & = \big\{\frac{1}{n}\sum_{i = 1}^n (y_i-\psi)(y_i-\psi)^{\T}\big\}- t^2(\bar{y}-\psi)(\bar{y}-\psi)^{\T} = A - t^2vv^{\T}.
\end{align}
These maximum likelihood estimators are valid as long as $\hat{\Lambda}^{-1}(t)$ is positive definite. The largest $t$ such that $\hat{\Lambda}^{-1}(t)$ is positive definite, $t_{max}$, is found by solving for $t$ in the equation $\vert \hat{\Lambda}^{-1}(t) \vert = 0$. Using the matrix determinant lemma we find that 
\begin{align}
\vert \hat{\Lambda}^{-1}(t) \vert = \det(A)(1 - t^2v^{\T}A^{-1}v).
\end{align}
The components of \eqref{densline} needed to compute the directional $p$-value are  
\begin{align}
\exp\big[\ell\{\hat{\varphi}(0);s(t)\} - \ell\{\hat{\varphi}(t);s(t)\}\big] &= \vert \hat{\Lambda}^{-1}(t) \vert^{n/2} , \\
|J_{\varphi\varphi}\{\hat{\varphi}(t);s(t)\}|^{-1/2} &= |\hat\Lambda^{-1}(t)|^{-(q+1)/2}.
\end{align}

As the canonical parameter is linear in $\Lambda$, the nuisance parameter information term does not depend on $t$ and can be ignored. The directional $p$-value is 
\begin{align}
\frac{\int^{(v^{\T}A^{-1}v)^{-1/2}}_1t^{p-1}(1 - t^2v^{\T}A^{-1}v)^{\frac{n - p - 2}{2}} \, dt}{\int^{(v^{\T}A^{-1}v)^{-1/2}}_0t^{p-1}(1 - t^2v^{\T}A^{-1}v)^{\frac{n - p - 2}{2}} \, dt}.
\end{align}

There is a striking similarity between this and the directional $p$-value in \S 3.1, and the same change of variables  can be applied here. By the Sherman-Morrison formula  
\begin{align}
\frac{1 - v^{\T}A^{-1}v}{v^{\T}A^{-1}v} = (\bar{y} - \psi)^{\T}\big\{\frac{1}{n}\sum_{i = 1}^n(y_i - \bar{y})(y_i - \bar{y})^{\T} \big\}^{-1}(\bar{y} - \psi) = v^{\T}B^{-1}v.
\label{partofhotstat}
\end{align}
After the appropriate change of variables \eqref{partofhotstat} appears in the bound of the integrals in the directional $p$-value. Eventually, we get that
\begin{align}
p(\psi) = 1 - G\big(\frac{n - p}{p} v^{\T} B^{-1}v\big),
\end{align}

where $G$ is the CDF of a $F(p, n - p)$ random variable. As Hotelling's $T^2$ statistic is equal to $(n - 1)v^{\T}B^{-1}v$ and  $\frac{n - p}{p (n - 1)}T^2 \sim F(p,n-p)$, this shows that the directional test is identical to Hotelling's $T^2$ test.

\section{Discussion}
\label{sec:disc}

Directional testing is an approach to testing of a multivariate parameter that in effect creates a one-dimensional sub-model, by restricting attention to the line on the parameter space that is dictated by the observed data. The computation of these tests has recently been simplified by relying on saddlepoint approximations to the distributions, rather than computing them exactly.
\\
\\
This work concentrates on models for which saddlepoint approximations are exact or nearly exact, and shows that conventional $F$-tests emerge from the directional approach. This helps to explain the accuracy of the tests demonstrated in \citet{Davison.etal:2014} and \citet{Fraser.etal:2016}. 
\\
\\
It is interesting to note that all of the directional $p$-value integrands appearing in this work share a common structure. Let  $\hat{\sigma}^2(t)$ be a measure of variability of observing $s(t)$ under $H_{\psi}$.  The integrand of the directional $p$-values are given by $t^{d - 1}\hat{\sigma}^2(t)^{\alpha/ 2}$ where $\alpha$ depends on $n$,$d$ and $p$. Small directional $p$-values correspond to observed data that have a relatively high weighted variability estimate under $H_{\psi}$.
\\
\\
The hypotheses considered in \S 3 constrain the mean vector to a linear subspace of the parameter space, and are also invariant under affine transformations of the parameter. The $F$-tests in \S 3 are derived as   most powerful invariant tests in \citet[Ch.7]{Lehmann.Romano:2005}, and shown there to effectively test a scalar parameter, the noncentrality parameter of the related $F$ distribution. 
\\
\\
In ongoing work we are preparing an R package to construct directional tests in generalized linear models, including gamma, Poisson, and logistic regression. With discrete probability functions, such as the binomial and Poisson, saddlepoint methods are not exact because they are implicitly continous, but simulation results in Example 4.2 in \citet{Davison.etal:2014} and in work in progress indicates that the directional $p$-values continue to be very accurate. 
\\
\\
It is also straightforward to develop directional tests for normal theory non-linear regression models of the form $y_i \sim N\{\eta_i(\beta),\sigma^2\}, i = 1, \dots, n$, but we have been unable to verify that these are the same as the conventional $F$-test based on the tangent model approximation to the mean surface.
\\
\\
By their nature, we might expect directional tests to have low power in regions of the parameter space that are not suggested by the data. This point was raised in work as yet unpublished by Jensen. However for at least some settings the work here shows that the tests are the same as conventional $F$-tests for multivariate hypotheses, so share their power properties.

\section*{Acknowledgements}
This research was supported in part by the Natural Sciences and Engineering Research Council of Canada and by the Italian Ministry of Education under the PRIN 2015 grant 2015EASZF\_003 and the University of Padova (PRAT 2015 CPDA153257). We would like to thank D.A.S. Fraser and I. Kosmidis for helpful discussion.

\newpage
\begin{center}
{\large\bf SUPPORTING INFORMATION}
\end{center}
Additional information about the saddlepoint approximation and directional tests, as well as detailed calculations to support the analytical results in \S 2 and \S 3 are provided in the supplementary material at the end of this document.

\newpage

\begin{center}
\LARGE Supporting Information for \\ \bf A New Look at $F$-tests
\end{center}

\section*{S.1 Introduction}

The tests examined in the main paper use a directional argument proposed in \citet{Fraser.Massam:1985} and developed further, using higher order approximation theory, in \citet{Davison.etal:2014} and \citet{Fraser.etal:2016}. For completeness we provide a summary of directional tests in \S S.1, and provide the  formulae needed for the examples in \S S.2 - S.4. 

\section*{S.2 Directional testing}
\subsection*{S.2.1 A model on $\mathbb{R}^p$}
Suppose our model for $y = (y_1, \dots, y_n)$ is a linear exponential family
\begin{equation}
\label{eq1}
f(y;\theta) = \exp[\varphi^\T(\theta)u(y)-\kappa\{\varphi(\theta)\}]d(y) ,\quad\varphi\in\mathbb{R}^p
\end{equation}
with sufficient statistic $u$ and canonical parameter $\varphi$. Inference for $\varphi$ is based on the marginal distribution of $u$, which is again an exponential family
\begin{equation}\label{eq2}
f(u;\theta) = \exp[\varphi^\T(\theta)u - \kappa\{\varphi(\theta)\}]\tilde d(u).
\end{equation}
The function $\tilde d(\cdot)$ is obtained by marginalizing (\ref{eq1}) and may not be available explicitly, but the saddlepoint approximation to the density of $u$ has relative error $O(n^{-3/2})$ in continuous models:
\begin{equation}\label{eq3}
f_{SP}(u;\theta) = \frac{e^{k/n}}{(2\pi)^{p/2}} |\hat\jmath|^{-1/2}\exp[\ell\{\varphi(\theta);u\}-\ell\{\varphi(\hat\theta);u\}],
\end{equation}
where $\ell(\varphi;u) = \varphi^{\T}u - \kappa(\varphi)$ is the log-likelihood function, $\hat\jmath = -\partial^2\ell(\hat\varphi)/\partial\varphi\partial\varphi^{\T}$ is the observed Fisher information, $\hat\theta$ is the maximum likelihood estimator, and $\exp(k/n)/(2\pi)^{p/2}$ is an approximation to the normalizing constant. 

Directional tests take as their starting point the approximation (\ref{eq3}). If the originating model is not in the exponential family, then an approximation to it, the tangent exponential model, is used instead. The construction of the tangent exponential model and its saddlepoint approximation are described in the Appendix of \citet[]{Fraser.etal:2016}; see in particular Eq.~(A2). Since the examples in the current paper are all exponential family models, this step is not needed.

\subsection*{S.2.2 Nuisance parameters}

Suppose that the parameter $d$-dimensional parameter  $\theta$ in our original model can be partitioned as $\theta = (\psi, \lambda)$, where $\psi$ is a $d$-dimensional parameter of interest and $\lambda$ is a $p-d$-dimensional nuisance parameter. Denote by $\hat\lambda_\psi$ the constrained maximum likelihood estimator of $\lambda$ when $\psi$ is fixed. 

In the special case that $\theta = \varphi$, so that the parameter of interest is a sub-vector of the canonical parameter, and if the original model is a full exponential family, then 
\begin{equation}\label{eq4}
f(u_1 , u_2;\psi,\lambda) \propto \exp\{\psi^{\T}u_1 + \lambda^{\T}u_2-\kappa(\psi,\lambda)\}\tilde d(u),
\end{equation}
and the conditional distribution of $u_1$ given $u_2$ is free of $\lambda$. The saddlepoint approximation can be used again to give an accurate approximation to the conditional density. Directional testing for models of this form are developed and illustrated in \citet{Davison.etal:2014}. 

If the parameter of interest is not a linear function of the canonical parameter, in which case we write $\psi = \psi(\varphi)$, such a reduction by conditioning is not available. None-the-less, it can be verified that there is a unique variable that measures $\psi$, and that this variable is obtained by constraining the sufficient statistic $u$ to the $d$-dimensional sample space obtained by fixing the constrained maximum likelihood estimate of the nuisance parameter to its observed value. The saddlepoint approximation to the density of this variable is 
\begin{equation}\label{eq5}
h(s;\psi) = \frac{\exp(k'/n)}{(2\pi)^{d/2}}\exp\{\ell(\hat\varphi_\psi;s)-\ell(\hat\varphi(s)\}|\hat J_{\varphi\varphi}|^{-1/2} |\tilde J_{(\lambda\lambda)}|^{1/2}, \quad s \in L_\psi, 
\end{equation}
where $L_\psi$ is the plane in the sample space with $\hat\lambda_\psi$ fixed, so that $h$ above is a density on $\mathbb{R}^d$. 
In (\ref{eq5}), $s = u - u^0$ is a centred version of the sufficient statistic, and $\ell(\varphi;s) = \varphi^{\T}s + \ell^0(\varphi)$ is an exponential tilt of the observed log-likelihood function $\ell(\theta;y^0)$ 
in the original model. The centering is described in detail in \citet[\S3.1]{Davison.etal:2014} and assumed in \citet{Fraser.etal:2016}.
The determinants in (\ref{eq5}) are:
\begin{equation}\label{Fisherinfo}
|\hat J_{\varphi\varphi}| = |J_{\varphi\varphi}\{\hat\varphi(s)\}| = |-\left. \partial^2\ell(\varphi;s)/(\partial\varphi\partial\varphi^{\T})\right |_{\varphi = \hat\varphi(s)}, 
\end{equation}
and 
\begin{equation}\label{Fisherinfo2}
|\tilde J_{(\lambda\lambda)}| = |J_{(\lambda\lambda)}(\hat\varphi_\psi)| =  \left |-\frac{\partial^2 \ell(\hat\varphi_\psi;s)}{\partial\lambda\partial\lambda^{\T} }\right| \left |\frac{\partial\varphi(\hat\theta_\psi)}{\partial\lambda }\right |^{-2}.
\end{equation}
The second determinant is not needed when the parameter of interest is linear in $\varphi$. However, it turns out to be independent of $t$ in the examples in \S 2.1, 2.3 and 3.2, even though the parameter of interest is not linear in the canonical parameter. In these cases the canonical parameter is a linear function of the nuisance parameter and so \eqref{Fisherinfo2} does not depend on $t$.

\subsection*{S.2.3 Directional testing}
The directional test of the hypothesis $\psi(\varphi) = \psi$ is carried out in $L_\psi$ by finding the line that joins $s^0$ with the value of $s$, call it $s_\psi$, that would give $\hat\varphi_\psi$ as the maximum likelihood estimate of the parameter. The observed value $s^0$ gives $\hat\varphi$ as the maximum likelihood estimate. The $p$-value is computed as the probability of $s$ being larger than the observed value $s_0$, on the line between the two values $s_{\psi}$ and $s^0$. Another way to describe it is that we measure the magnitude of the vector $s_0 - s_\psi$, in $L_\psi$, conditional on its direction. This gives a one-dimensional measure of how much ``larger'' the observed value $s^0$ is than would be expected under the hypothesis. We parameterize this line in the sample space by $t$, and because we center the sufficient statistic so that $s^0=0$, the line is simply $s(t) = s_\psi + t(s^0-s_\psi) = (1-t)s_\psi $. The directional $p$-value is then 
\begin{equation}\label{eq6}
p(\psi) = \frac {\int_1^{\infty} t^{d-1}h(t;\psi)dt}{\int_0^{\infty} t^{d-1}h(t;\psi)dt},
\end{equation}
where $h(t;\psi) = h\{s(t);\psi\}$ using (\ref{eq4}), and the inflation factor $t^{d-1}$ comes from the Jacobian of the transformation to polar coordinates.

\section*{S.3 Ratio of exponential rates}

\subsection*{S.3.1 Finding the integrand of the directional $p$-value}
We consider the directional test for the null hypothesis $H_{\psi}$ that $\theta_1/\theta_2 = \psi$ where $y_{ij} \sim \text{exp}(\theta_j)$, $j = 1,2 \; i = 1,...,n_i$. Note that this test is slightly different than the test performed Davison et al. (2014) as here we are testing the ratio of rates. Under $H_{\psi}$ the constrained MLE is 
\begin{align*}
\hat{\theta}_{1\psi} & = \frac{n\psi}{u_1\psi + u_2},
\\
\hat{\theta}_{2\psi} & = \frac{n}{u_1\psi + u_2},
\end{align*}
where $n = n_1 + n_2$ and $u_j = \Sigma_i y_{ij}$. By solving the score equation it is found that the value of $u_j$ that has $\hat{\theta}_{\psi}$ as its global MLE is $u_{j\psi} = n_j / \hat{\theta}_{j\psi}$. The line between $u_{\psi}$ and the observed value of $u$, $u^0 = (u_1^0,u_2^0)$ is
\begin{align*}
u_1(t) = & \frac{n_1}{n}(u_1^0 + \frac{u_2^0}{\psi})+t\big\{u_1^0 - \frac{n_1}{n}(u_1^0 + \frac{u_2^0}{\psi}) \big\} = u_{1\psi} + t(u_1^0 - u_{1\psi}),
\\
u_2(t) = & \frac{n_2}{n}(u_1^0\psi + u_2^0) + t\big\{u_2^0 -  \frac{n_2}{n}(u_1^0\psi + u_2^0)\big\} = u_{2\psi} + t(u_2^0 - u_{2\psi}),
\end{align*}
so that 
\begin{align*}
\exp[\ell\{\hat{\varphi}_{\psi};s(t)\} - \ell\{\hat{\varphi};s(t)\}] \propto \exp\bigg\{-\frac{n_1u_1(t)}{u_{1\psi}} - \frac{n_2u_2(t)}{u_{2\psi}}\bigg\}u_1(t)^{n_1}u_2(t)^{n_2}.
\end{align*}
The determinant of the Hessian of the negative log-likelihood with respect to $\varphi = (\theta_1, \theta_2)$ is the determinant of the Hessian of $-n_1\log(\theta_1) - n_2\log(\theta_2)$. This is easily found to be $\vert J_{\varphi\varphi}\vert = n_1 n_2 /(\theta_1\theta_2)^2$. We have that $\hat{\theta}_i(t) = n_i / u_i(t)$ and thus 
\begin{align*}
\vert J_{\varphi\varphi} \vert^{-\frac{1}{2}}\exp[\ell\{\hat{\varphi}_{\psi};s(t)\} - \ell\{\hat{\varphi};s(t)\}] \propto \exp\bigg\{-\frac{n_1u_1(t)}{u_{1\psi}} - \frac{n_2u_2(t)}{u_{2\psi}}\bigg\}u_1(t)^{n_1-1}u_2(t)^{n_2-1}.
\end{align*}
Next we check the nuisance parameter adjustment term. Formulating the log-likelihood in terms of nuisance parameter $\theta_2 = \lambda$ we get
\begin{align*}
\ell(\psi,\lambda) = -u_1(t)\psi\lambda - u_2(t)\lambda + \ell^0(\psi,\lambda).
\end{align*}
It is clear that the second derivative of this function with respect to $\lambda$ only contains terms that do not involve $t$. The dimension of our parameter of interest in this case is $d = 1$ so that $t^{d-1} = 1$. The sufficient statistic $u(t)$ is viable as long as it remains positive. As a result, $t_{max}$ is the largest $t$ such that both $u_1(t)$ and $u_2(t)$ are non-negative. Notice that  $u_1^0 - u_{1\psi} > 0$ if and only if $\psi\bar{y}_1 \geq \bar{y}_2$ and likewise $u_2^0 - u_{2\psi} > 0$ if and only if $\psi\bar{y}_1 \leq \bar{y}_2$. so that 
\begin{align*}
t_{max} = \frac{u_{1\psi}}{u_{1\psi} - u_1^0}\text{I}(\psi\bar{y}_1 \leq \bar{y}_2) + \frac{u_{2\psi}}{u_{2\psi} - u_2^0}\text{I}(\psi\bar{y}_1 \geq \bar{y}_2).
\end{align*}
We let $a_j$ be equal to the quantity $u_{j\psi}/(u_{j\psi} - u_j^0)$. It can be shown that $(u_{1\psi} - u_1^0)/(u_{2\psi} - u_2^0) = -1/\psi$ and $u_{2\psi}/u_{1\psi} = \psi n_2/ n_1$, which shows that $n_1a_2 +n_2a_1 = 0$. Now 
\begin{align*}
\exp\bigg\{-\frac{n_1u_1(t)}{u_{1\psi}} - \frac{n_2u_2(t)}{u_{2\psi}}\bigg\} \propto \exp\bigg\{ t\big(\frac{n_1}{a_1} +  \frac{n_2}{a_2}\big)\bigg\} = 1,
\end{align*}
and thus the directional $p$-value is given by
\begin{align}
p(\psi) = \frac{\int_{1}^{t_{max}}u_1(t)^{n_1-1}u_2(t)^{n_2-1}dt}{\int_0^{t_{max}}u_1(t)^{n_1-1}u_2(t)^{n_2-1}dt} = \frac{\int_{1}^{t_{max}}(1 - \frac{t}{a_1})^{n_1-1}(1-\frac{t}{a_2})^{n_2-1}dt}{\int_0^{t_{max}}(1 - \frac{t}{a_1})^{n_1-1}(1-\frac{t}{a_2})^{n_2-1}dt}.
\label{dirpval}
\end{align}

\subsection*{S.3.2 Making a change of variables}
Assume that $t_{max} = a_1$, so that $\psi\bar{y}_1/\bar{y}_2 \leq 1$. The numerator of \eqref{dirpval} can be written as 
\begin{align*}
p_{num} = & \int^{a_1}_1 (1 - \frac{t}{a_1})^{n_1 + n_2} \bigg(\frac{1 - \frac{t}{a_2}}{1 - \frac{t}{a_1}}\bigg)^{n_2 - 1} (1 - \frac{t}{a_1})^{-2} dt.
\end{align*}
Make the change of variables $x = (1 - t/a_2)/(1- t/a_1)$, we then get that
\begin{align*}
p_{num} = & k_1\int_{\frac{1 - 1/a_2}{1 - 1 / a_1}}^{\infty} x^{n_2 - 1}\bigg( \frac{a_1 - a_2}{a_1 - a_2 x} \bigg)^{-n_1 - n_2}dx = k_2\int_{\frac{1 - 1/a_2}{1 - 1 / a_1}}^{\infty} x^{n_2 - 1}\big( 1 - \frac{a_2x}{a_1} \big)^{-n_1 - n_2}dx.
\end{align*}
We know from a previous calculation that $a_2/a_1 = -n_2/n_1$. Furthermore,
\begin{align*}
\frac{1 - 1/a_1}{1-1/a_2} = \frac{(a_2 - 1)a_1}{(a_1 - 1)a_2} = -\frac{u_2^0(u_{1\psi} - u_1^0)n_1}{u_1^0(u_{2\psi} - u_2^0)n_2} = \frac{\bar{y}_2}{\bar{y}_1\psi}.
\end{align*}
Therefore we see that
\begin{align}
p_{num} = k_2 \int_{\frac{\bar{y_2}}{\bar{y_1}\psi}}^{\infty} x^{n_2 - 1}\big( 1 + \frac{n_2}{n_1}x \big)^{-n_1 - n_2}dx.
\label{pvalintegrand}
\end{align}
Now if we perform the same change of variables on the integral in the denominator of the directional $p$-value the same normalizing constant $k_2$ will be produced. The bounds of the integral will change as follows
\begin{align*}
\int^{a_1}_0 \Rightarrow \int_1^{\infty}.
\end{align*}
The integrand in \eqref{pvalintegrand} is the density of a $F(2n_2,2n_1)$ random variable up to a normalizing constant. Thus if $W \sim F(2n_2,2n_1)$
\begin{align*}
p(\psi) = \frac{\text{P}_{W}(W > \frac{\bar{y}_2}{\bar{y}_1\psi})}{\text{P}_{W}(W > 1)}.
\end{align*}
Similarly when $\bar{y}_1\psi/\bar{y}_2 \geq 1$ we get 
\begin{align*}
p(\psi) = \frac{\text{P}_{W}(W< \frac{\bar{y}_1\psi}{\bar{y}_2})}{\text{P}_{W}(W < 1)}.
\end{align*}
Consequently, the directional test is identical to that of the appropriate F-test. In summary, the directional $p$-value is
\begin{align*}
p(\psi) = \text{I}\big(\bar{y}_1\psi < \bar{y}_2\big)\frac{\text{P}_{W}(W > \frac{\bar{y}_2}{\bar{y}_1\psi})}{\text{P}_{W}(W > 1)} + \text{I}\big(\bar{y}_1\psi > \bar{y}_2 \big) \frac{\text{P}_{W}(W < \frac{\bar{y}_1\psi}{\bar{y}_2})}{\text{P}_{W}(W < 1)}.
\end{align*}

\section*{S.4 Ratio of normal variances}
Suppose that $y_{ij}\sim N(\mu_i, \sigma_i^2)$ are independent random variables for $i=1,2$ and $j=1\dotsc,n_i$, and we wish to test $H_{\psi}: \sigma_1^2/\sigma_2^2 = \psi$.
A computation very similar to that given in the previous section shows that the integrand of the directional $p$-value for testing $H_{\psi}$ is 
\begin{align}
(1-tb_1)^{(n_1-3)/2}(1-tb_2)^{(n_2-3)/2},
\label{varintegrand}
\end{align}
with $b_i=(\hat{\sigma}_{i\psi}^2 -v_i^2)/\hat{\sigma}_{i\psi}^2$. The biased within-group sample variances are $v_i^2$ while $\hat{\sigma}_{i\psi}^2$ is the constrained maximum likelihood estimator for $\sigma_i^2$. There is a clear resemblance between the integrands of \eqref{dirpval} and \eqref{varintegrand}. The same change of variables used for the exponential rates example can be used here with the only minor difference being that we set $a_i = 1/b_i$ so that \eqref{varintegrand} equals $(1-t/a_1)^{(n_1-3)/2}(1-t/a_2)^{(n_2-3)/2}$. All the bounds of the integrals will change in the same way as in the previous example. In particular, we find that
\begin{align*}
\frac{1 - 1/a_2}{1 - 1/a_1} = \frac{1 - b_2}{1 - b_1} = \frac{\hat{\sigma}_{1\psi}^2 v_2^2}{\hat{\sigma}_{2\psi}^2 v_1^2} = \frac{\psi v_2^2}{ v_1^2}.
\end{align*}
If $v_1^2 \leq \psi v_2^2$ so that $t_{max} = 1/a_1$ we get that 
\begin{align*}
p_{num} = \int^\infty_{\psi v_2^2/v_1^2} x^{\frac{n_2 - 1}{2} - 1}\big(1 + \frac{n_2}{n_1}x\big)^{-\frac{(n_1 - 1) + (n_2 - 1)}{2}} \, dt.
\end{align*}
This integrand is not quite the density of a F-distribution due to the factor of $n_2/n_1$ appearing instead of $(n_2 - 1)/ (n_1 - 1)$. To fix this we make the additional change of variables in both $p_{num}$ and $p_{den}$ from $x$ to \\ $\{n_1(n_2 - 1)\}/\{n_2(n_1 - 1)\}x$. If $W \sim F(n_2 - 1, n_1 - 1)$ and $s_i^2$ are the unbiased sample variances we get the desired result that 
\begin{align*}
p(\psi) = \frac{P_W\big(W > \frac{\psi s_2^2}{s_1^2}\big)}{P_W\big(W > \frac{n_2(n_1 - 1)}{n_1(n_2 - 1)}\big)}.
\end{align*}
The case where $v_1^2 > \psi v_2^2$ is handled similarly. 

\section* {S.5 Linear regression with linear constraints}

\subsection* {S.5.1 Calculating $\exp [\ell\{ \hat{\varphi}(0);s(t)\} - \ell\{ \hat{\varphi}(t);s(t)\}]$}

 Let $y_{i} \sim N(x_i^\T\beta, \sigma^2)$, $ i = 1,\dots,n$, where $y_{i}^0$ are realizations of $y_{i}$ and all of the $y_i$'s are independent. Both $x_i$ and $\beta$ are $p \times 1$ vectors and $\sigma^2$ is an unknown nuisance parameter. We form the $n \times p$ matrix $X$ by taking i'th row of $X$ to be $x_i$. Here we wish to test $H_{\psi}: A \beta = \psi$ using the directional test. The matrix $A$ has dimension $d \times p$ and is of rank $d$ which ensures that the linear constraint is not redundant. The constrained maximum likelihood estimator for $\beta$ under $H_{\psi}$ can be found using Lagrange multipliers and is given by
\begin{align*}
\hat{\beta}_{\psi} & = \hat{\beta} - (X^\T X)^{-1}A^\T\{A(X^\T X)^{-1}A^\T\}^{-1}(A\hat{\beta} - \psi) 
\\
& = \hat{\beta} - \frac{1}{2}(X^\T X)^{-1}A^\T\hat{\lambda}.
\end{align*}
The Lagrange multiplier equation used to find this constrained MLE also yields
\begin{align*}
\hat{\beta}_{\psi}^\T X^\T(y - X\hat{\beta}_{\psi}) & = \frac{1}{2}\psi^\T\hat{\lambda}
\\
& = \psi^\T\{A(X^\T X)^{-1}A^\T\}^{-1}(A\hat{\beta} - \psi).
\end{align*}
The constrained MLE for $\sigma^2$ is just the average sum of squared error under $\hat{\beta}_{\psi}$. 
The log-likelihood in this situation is 
\begin{align*}
\ell(\beta,\sigma^2) & =  -\frac{y^\T y}{2\sigma^2} + \frac{y^\T X\beta}{\sigma^2} - \frac{\beta^\T X^\T X\beta}{2\sigma^2} - \frac{n}{2}\log(\sigma^2)
\\
& = \begin{bmatrix}
y^\T y & y^\T X
\end{bmatrix}
\begin{bmatrix}
-\frac{1}{2\sigma^2}
\\
\frac{\beta}{\sigma^2}
\end{bmatrix}
 - \kappa(\beta,\sigma^2)
 \\
 & = u^\T(y)
\varphi - \kappa(\varphi).
\end{align*}
The sufficient statistics here are $y^\T y$ and $y^\T X$. These sufficient statistics have unconstrained MLEs that are equal to the constrained MLE when they solve the following equations:
\begin{align*}
(X^\T X)\hat{\beta}_{\psi} = X^\T y
\\
\text{and}
\\
\frac{1}{n}(y - X\bcon)^\T(y- X\bcon) = \hat{\sigma}^2_{\psi}
\\
\implies \frac{1}{n}(y^\T y - \bcon^\T X^\T X\bcon) = \hat{\sigma}^2_{\psi}
\\
y^\T y = n \hat{\sigma}^2_{\psi} + \bcon^\T X^\T X\bcon.
\end{align*}
Thus 
\begin{align*}
s_{\psi} = \begin{bmatrix}
n \hat{\sigma}^2_{\psi} + \bcon^\T X^\T X\bcon -    (y^0)^{\T}y^0
\\
(X^\T X)\bcon - X^\T y^0
\end{bmatrix}.
\end{align*}
We define $s(t) = (1-t) s_{\psi}$. Then let $u(t) = u^0 + s(t)$ so that $\ell(\varphi;t) = u^\T(t)\varphi - \kappa(\varphi)$. We let $u_1(t)$ be the first entry of $u(t)$ and $u_2(t)$ be the remaining entries of $u(t)$. We see that
\begin{align*}
\frac{\partial}{\partial\sigma^2}\ell(\varphi;t) & = \frac{1}{\sigma^4}\big\{\frac{u_1(t)}{2} - u_2(t)\beta + \frac{\beta^\T X^\T X\beta}{2}\big\} - \frac{n}{2\sigma^2} = 0
\\
\implies \hat{\sigma}^2(t) & = \frac{1}{n}\big\{u_1(t) - 2u_2(t)\hat{\beta}(t) + \hat{\beta}(t)^\T X^\T X\hat{\beta}(t)\big\}.
\\
\frac{\partial}{\partial \beta}\ell(\varphi;t) & = \frac{1}{\sigma^2}\big\{u_2(t) - X^\T X\beta\big\} = 0
\\
\implies \hat{\beta}(t) & = (X^\T X)^{-1}u_2(t).
\end{align*}
We now find formulas for the log-likelihood terms appearing in the exponent density used for the directional $p$-value calculation. As the directional $p$-value takes a ratio of such densities we can ignore all multiplicative factors not involving $t$ in the subsequent calculations:
\begin{align*}
\ell\{\hat{\varphi}(0);s(t)\} & \propto \frac{1}{\hat{\sigma}^2(0)}\big\{-\frac{1}{2}u_1(t) + u_2^\T(t)\hat{\beta}(0)\big\}
\\
& \propto \frac{1}{\hat{\sigma}^2(0)}\bigg[\frac{t}{2}\big\{n \hat{\sigma}^2_{\psi} + \bcon^\T X^\T X\bcon -    (y^0)^{\T}y^0\big\} + t\big\{(y^0)^{\T}X - \bcon X^\T X\big\}\bcon \bigg]
\\
& \propto  \frac{1}{\hat{\sigma}^2(0)}\bigg(\frac{nt}{2}\big[ \hat{\sigma}^2_{\psi}- \{(y^0)^{\T}y^0 - 2(y^0)^{\T}X\bcon + \bcon^\T X^\T X\bcon\} \big] \bigg) = 0.
\end{align*}
Similarly,
\begin{align*}
\ell\{\hat{\varphi}(t);s(t)\} & = \frac{1}{\hat{\sigma}^2(t)}\big\{-\frac{1}{2}u_1(t) + u_2^\T(t)\hat{\beta}(t)\big\} - \frac{\hat{\beta}^\T(t) X^\T X\hat{\beta}(t)}{2\hat{\sigma}^2(t)} - \frac{n}{2}\log\{\hat{\sigma}^2(t)\}
\\
& = \frac{1}{\hat{\sigma}^2(t)}\big\{-\frac{1}{2}u_1(t) + \frac{1}{2}u_2^\T(t)\hat{\beta}(t)\big\} - \frac{n}{2}\log\{\hat{\sigma}^2(t)\}
\\
& = -1 - \frac{n}{2}\log\{\hat{\sigma}^2(t)\}.
\end{align*}
Consequently, $\exp\big[\ell\{\hat{\varphi}(0);s(t)\} - \ell\{\hat{\varphi}(t);s(t)\}\big] = \{\hat{\sigma}^2(t)\}^{\frac{n}{2}}$.

\subsection*{S.5.2 Finding $\vert J_{\varphi\varphi}\{\hat{\varphi}(t) ; s(t)\}\vert$ and the nuisance parameter adjustment}

To start we find all of the second order derivatives of $\kappa(\varphi)$ yielding

\begin{align*}
\frac{\partial^2\kappa}{\partial\varphi_{1}^2} & = \frac{n}{2\varphi_{1}^2} + \frac{1}{\varphi_{1}^3}\sumonn (\sum_{j = 2}^{p+1} X_{ij}\varphi_j)^2,
\\
\frac{\partial^2\kappa}{\partial\varphi_{1}\partial\varphi_k} & = -\frac{1}{\varphi_{1}^2}\sumonn X_{ik}\sum_{j = 2}^{p+1} X_{ij}\varphi_j,
\\
\frac{\partial^2\kappa}{\partial\varphi_{k}\partial\varphi_l} & = \frac{1}{\varphi_{1}}\sumonn X_{ik}X_{il}. 
\end{align*}

\noindent Define $\bar{\varphi}$ to be the vector containing the last $p$ entries of $\varphi$. From the second order derivatives above we find that the Hessian of the negative log-likelihood function has the form

\begin{align*}
J_{\varphi\varphi}\{\varphi;s(t)\} = \frac{1}{\varphi_{1}}\begin{bmatrix}
(\frac{n}{2\varphi_{1}} + \frac{1}{\varphi_{1}^2}\bar{\varphi}^\T X^\T X\bar{\varphi}) &  -\frac{1}{\varphi_{1}}\bar{\varphi}^\T X^\T X
\\
-\frac{1}{\varphi_{1}}X^\T X\bar{\varphi} & X^\T X
\end{bmatrix}.
\end{align*}

\noindent We multiply the above matrix on the left by the matrix

\begin{align*}
\begin{bmatrix}
1 & 0^\T
\\
0 & (X^\T X)^{-1}
\end{bmatrix}.
\end{align*}

\noindent In doing this the determinant of the resulting matrix only changes by the constant $\det\{(X^\T X)^{-1}\}$. We then find that find that:

\begin{align*}
\vert J_{\varphi\varphi}\{\varphi;s(t)\} \vert \propto (\frac{1}{\varphi_{1}})^{p+1}\det\bigg(\begin{bmatrix}
(\frac{n}{2\varphi_{1}} + \frac{1}{\varphi_{1}^2}\bar{\varphi}^\T X^\T X\bar{\varphi}) & -\frac{1}{\varphi_{1}}\bar{\varphi}^\T X^\T X
\\
-\frac{1}{\varphi_{1}}\bar{\varphi} & I_p
\end{bmatrix}\bigg).
\end{align*}

\noindent The determinant above can be found by performing a cofactor expansion along the first column of the matrix and then performing a cofactor expansion along first row of the resulting minor. Fortunately, performing this cofactor expansion twice on the i'th entry of the first column and the j'th entry of the first row will produce a minor that has a row of zeros if $i \neq j$. As a result, we only have to be concerned about when $ i = j$, but this case is simple as it is just minus one times the i'th entry of the last column times the i'th entry of the last row multiplied by the determinant of the identity. The bottom right entry has to be treated separately, but it clearly just returns itself in the cofactor expansion. In short

\begin{align*}
\vert J\{\varphi;s(t) \} \vert \propto & (\frac{1}{\varphi_{1}})^{p+1}\bigg\{-\frac{1}{\varphi_{1}^2}\bar{\varphi}^\T X^\T X\bar{\varphi} + (\frac{n}{2\varphi_{1}} + \frac{1}{\varphi_{1}^2}\bar{\varphi}^\T X^\T X \bar{\varphi})\bigg\} 
\\
\propto & (\frac{1}{\varphi_{1}})^{p+2}.
\end{align*}
At this point we are able to construct the integrand for the directional test:
\begin{align*}
t^{d-1}\vert J\{\hat{\varphi}(t);s(t)\}\vert^{-\frac{1}{2}}\exp{\big[\ell\{\hat{\varphi}(0) ; s(t)\} - \ell\{\hat{\varphi}(t) ; s(t) \}\big]} 
& \propto t^{d - 1} \hat{\sigma}^2(t)^{\frac{n-p-2}{2}}.
\end{align*}
This hypothesis is a linear hypothesis, meaning that the canonical parameter can be partitioned into the parameter of interest and a nuisance parameter. To see this, we note that $A\beta = \psi$ is equivalent to $A\beta/ \sigma^2 = \psi/ \sigma^2$. Consequently, $H_{\psi}$ holds if and only if
\begin{align}
\begin{bmatrix}
A & 2\psi
\end{bmatrix}
\varphi(\theta) = M\varphi(\theta) = 0.
\label{linearhypoth}
\end{align}
The dimension of $M$ is $d \times p + 1$. We define $\tilde{M}$ to be a matrix formed by adding $p + 1 - d$ rows to $M$ in a manner so that all of the rows of $\tilde{M}$ are linearly independent. We then see that $u^\T(y)\varphi(\theta) = u^\T(y) \tilde{M}^{-1}\tilde{M}\varphi(\theta)$ since $\tilde{M}$ is invertible, and we can redefine our canonical parameter to be $\tilde{M}\varphi(\theta)$. By \eqref{linearhypoth} this new canonical parameter can be partitioned into the parameter of interest and a nuisance parameter.  
Thus we are in the scenario covered in \citet{Davison.etal:2014}. No nuisance parameter adjustment term is required in this case.

\subsection*{S.5.3  Making a change of variables}

As a reminder 
\begin{align*}
\hat{\sigma}^2(t) & = \frac{1}{n}\{u_1(t) - u_2^\T(t)\hat{\beta}(t)\}
\\
& = \frac{1}{n}\{u_1(t) - u_2^\T(t)(X^\T X)^{-1}u_2(t)\}.
\end{align*}
After some algebra the term involving $t$ disappears from the above expression and we are left with 
\begin{align*}
\hat{\sigma}^2(t) & = \{\sigma_{\psi}^2 - \frac{t^2}{n}(y^0  - X\bcon)^\T X(X^\T X)^{-1}X^\T(y^0 - X\bcon) \}
\\
& = (a - t^2b).
\end{align*}
Thus $t_{max} = (\frac{a}{b})^{\frac{1}{2}}$. We find  the integral in the numerator to be

\begin{align*}
\int^{\sqrt{\frac{a}{b}}}_1 t^{d-1}(a-bt^2)^{\frac{n-p-2}{2}}dt = k \int^{\sqrt{\frac{a}{b}}}_1 \big(2\frac{b}{a}t\big)\big(\frac{b}{a}t^2\big)^{\frac{d-2}{2}}\big(1-\frac{b}{a}t^2\big)^{\frac{n-p-2}{2}}dt.
\end{align*}
 
\noindent Make the change of variables $x = \frac{b}{a}t^2$:

\begin{align*}
& = k\int^1_{\frac{b}{a}} x^{\frac{d-2}{2}}(1-x)^{\frac{n-p-2}{2}}dx.
\end{align*}

\noindent Make the change of variables $x = \frac{1}{z}$:

\begin{align*}
& = k\int^1_{\frac{a}{b}}-\frac{1}{z^2} \big(\frac{1}{z}\big)^{\frac{d-2}{2}}\big(1-\frac{1}{z}\big)^{\frac{n-p-2}{2}}dz 
\\
& = k\int^{\frac{a}{b}}_1 \big(\frac{1}{z}\big)^{\frac{d+2}{2}}\big(1-\frac{1}{z}\big)^{\frac{n-p-2}{2}}dz 
\\
& = k\int^{\frac{a}{b}}_1 \big(\frac{1}{z}\big)^{\frac{n+d-p}{2}}(z-1)^{\frac{n-p-2}{2}}dz.
\end{align*}

\noindent Make the change of variables $t = z - 1$:

\begin{align*}
& = k\int^{\frac{a}{b}-1}_0 \big(\frac{1}{t+1}\big)^{\frac{n+d-p}{2}}t^{\frac{n-p-2}{2}}dt 
\\
& = k\int^{\frac{a - b}{b}}_0 \big(\frac{1}{t+1}\big)^{(\frac{n-p}{2}+\frac{d}{2})}t^{\frac{n-p}{2}-1}dt.
\end{align*}

\noindent Make the final change of variables $t = \frac{n-p}{d}x$:

\begin{align*}
& = k^\prime \int^{\frac{(a - b)d}{b(n-p)}}_0 \big(\frac{1}{\frac{n-p}{d}x+1}\big)^{(\frac{n-p}{2}+\frac{d}{2})}\big(\frac{n-p}{d}x\big)^{\frac{n-p}{2}-1}dx
\\
& = k^{\prime\prime}\int^{\frac{(a - b)d}{b(n-p)}}_0 \frac{\Gamma(\frac{n-p}{2} + \frac{d}{2})}{\Gamma(\frac{n-p}{2})\Gamma(\frac{d}{2})}\frac{n-p}{d-1} \big(\frac{1}{\frac{n-p}{d}x+1}\big)^{(\frac{n-p}{2}+\frac{d}{2})}\big(\frac{n-p}{d}x\big)^{\frac{n-p}{2}-1}dx
\\
& = k^{\prime\prime}P_W\bigg(W > \frac{\frac{b}{d}}{\frac{a - b}{n-p}}\bigg).
\end{align*}

\noindent where $W \sim F(d,n - p)$. The above sequence of changes of variables is equivalent to making the single change $x = \{(n - p)a - dbt^2\}/(dbt^2)$. Now performing the exact same sequence of changes of variables on the integral in the denominator will result in a similar expression, however the bounds of the integral over the F-distribution will be different. The bounds of the integral will change as follows:

\begin{align*}
\int_{0}^{\sqrt{\frac{a}{b}}} \Rightarrow \int_0^1 \Rightarrow \int_1^\infty \Rightarrow \int_0^\infty.
\end{align*}

\noindent Thus the integral in the denominator of the directional $p$-value equals

\begin{align*}
\int^{\sqrt{\frac{a}{b}}}_1 t^{q-1}(a-bt^2)^{\frac{n-p-2}{2}}dt = k^{\prime\prime} P_W(W > 0) = k^{\prime\prime}.
\end{align*}
We find $n(a-b)$ to be
\begin{align*}
(y^0)^{\T}y^0 - (y^0)^{\T}X(X^\T X)^{-1}X^\T y^0 + 2(X\hat{\beta} - y^0)^\T X\bcon.
\end{align*}
We know that $(y^0 - X\hat{\beta})^\T X\hat{\beta} = 0$. Thus
\begin{align*}
(X\hat{\beta} - y^0)^\T X\bcon & = \frac{1}{n}(y^0)^{\T}\big\{I_n - X(X^\T X)^{-1}X^\T\big\}X(X^\T X)^{-1}A^\T\{A(X^\T X)^{-1}A^\T\}^{-1}(A\hat{\beta} - \psi) 
\\
& = 0
\\
\implies (a-b) & = (y^0)^{\T}\{I_n - X(X^\T X)^{-1}X^\T\}y^0.
\end{align*}
This is simply the mean squared error under $\hat{\beta}$. Now we simplify $nb$ as
\begin{align*}
&(y^0  - X\bcon)^\T X(X^\T X)^{-1}X^\T (y^0 - X\bcon)  = (y^0)^{\T}X(X^\T X)^{-1}X^\T y^0 - 2(y^0)^{\T}X\bcon + \bcon^\T X^\T X\bcon
\\
 & = (y^0)^{\T}(X\hat{\beta} - X\bcon) - (y^0)^{\T}X\bcon + \bcon^\T X^\T X\bcon
 \\
 & = \frac{1}{2}(y^0)^{\T} X(X^\T X)^{-1}A^\T \hat{\lambda} - (y^0)^{\T} X\hat{\beta} + \frac{1}{2}(y^0)^{\T} X(X^\T X)^{-1}A^\T \hat{\lambda} + \bcon^\T X^\T X\bcon
\\
 & = (y^0)^{\T} X(X^\T X)^{-1}A^\T\hat{\lambda}- (y^0)^{\T} X\hat{\beta} + \hat{\beta}X^\T X\hat{\beta} - \hat{\beta}^\T X^\T X(X^\T X)^{-1}A^\T\hat{\lambda} + \frac{1}{4}\hat{\lambda}^\T A(X^\T X)^{-1}A^\T\hat{\lambda}
 \\
 & = \frac{1}{4}\hat{\lambda}^\T A(X^\T X)^{-1}A^\T\hat{\lambda}
 \\
 & = (A\hat{\beta} - \psi)^\T\big\{A(X^\T X)^{-1}A^\T\big\}^{-1}(A\hat{\beta} - \psi).
\end{align*}
Finally we find that the directional $p$-value equals the desired quantity:
\begin{align*}
p(\psi) = \text{P}_{W}\bigg\{W \geq \frac{(A\hat{\beta} - \psi)^\T\{A(X^\T X)^{-1}A^\T\}^{-1}(A\hat{\beta} - \psi)/d}{(y^0)^{\T} \{I - X(X^\T X)X\}y^0/(n-p)}\bigg\}.
\end{align*}

\section*{S.6 Multivariate normal mean}

\subsection*{S.6.1 Calculating $\exp [\ell\{ \hat{\varphi}(0);s(t)\} - \ell\{ \hat{\varphi}(t);s(t)\}]$}
Let $y_i \sim N_p(\mu, \Lambda^{-1})$ be $n$ observations from a multivariate Gaussian distribution with unknown concentration matrix $\Lambda$. Here we wish to test $H_{\psi}$ : $\mu = \psi$ using directional testing. The parameter of interest is $\mu$, while the nuisance parameter is $\con$. The log-likelihood of these observations is given by
\begin{align*}
\ell(\mu,\Lambda) & = -\frac{n}{2}\log(\vert \Lambda^{-1} \vert) -\frac{1}{2}\sumonn(y_i - \mu)^\T \con (y_i-\mu)
\\
& = \frac{n}{2}\log(\vert \con \vert) - \frac{1}{2}\sumonn y_i^\T \con y_i + \frac{1}{2}\mu^\T \con n\bar{y} + \frac{1}{2}n\bar{y}^\T \con \mu  -\frac{n}{2}\mu^\T \con \mu.
\end{align*}
By using the fact that the trace of a product of matrices is invariant under cyclic permutations and the trace of the product of two square matrices is the dot product of the vectorization of these matrices we can rewrite the log-likelihood as

\begin{align}
\ell(\mu, \con) = \begin{bmatrix}
\con\mu
\\
\text{vec}(\con)
\end{bmatrix}^\T \begin{bmatrix}
n \bar{y}
\\
\text{vec}(-\frac{1}{2}\sumonn y_i y_i^\T)
\end{bmatrix} + \frac{n}{2}\log(\vert \con \vert) - \frac{n}{2}\mu^\T \con\con^{-1}\con \mu
\label{likhoodpar}
\end{align}
We can rewrite the above log-likelihood in terms of the canonical parameter, $\varphi$, and sufficient statistic $u$:
\begin{align*}
\ell\{\varphi;u(y)\} = & \varphi^\T u + \frac{n}{2}\log(\vert \varphi_2 \vert) - \frac{n}{2}\varphi_1^\T \varphi_2^{-1} \varphi_1,
\\
\varphi = & \begin{bmatrix}
\varphi_1
\\
\varphi_2
\end{bmatrix}
= \begin{bmatrix}
\con \mu
\\
\con
\end{bmatrix}
= \begin{bmatrix}
\lambda\psi
\\
\lambda
\end{bmatrix},
\\
\theta = & \begin{bmatrix}
\psi
\\
\lambda
\end{bmatrix}
= \begin{bmatrix}
\mu
\\
\con
\end{bmatrix}.
\end{align*}
Throughout, we will be treat $\varphi_2$ as both a matrix and the vectorization of a matrix depending upon the context it is used in. The constrained MLE under $H_{\psi}$ is found by simply maximizing $\ell(\mu,\con)$ with respect to $\con$ while setting $\mu = \psi$. This yields the standard covariance matrix estimate $\hat{\con}_{\psi}^{-1} = \frac{1}{n} \sumonn (y_i - \psi)(y_i - \psi)^\T$. This matrix as well as $\psi$ will be used to find $s(t)$.
\\
\\
Finding $s(t)$ amounts to first finding a vector $s_{\psi}$ which when added to the observed value of the sufficient statistic has the constrained MLE, $\hat{\varphi}_{\psi}$, as its MLE. Once this is found $s(t)$ is given by $s(t) = (1-t)s_{\psi}$. The partial derivatives of $\ell(\varphi, s)$ with respect to $\varphi$ are
\begin{align*}
\frac{\partial \ell}{\partial \varphi_1} & = (u_1 + s_1) - n\varphi_2^{-1}\varphi_1
\\
\implies \qquad s_1 & = n\psi - u_1,
\\
\frac{\partial \ell}{\partial \varphi_2} & = (u_2 + s_2) + \frac{n}{2}\varphi_2^{-1} + \frac{n}{2}\varphi_2^{-1}\varphi_1\varphi_1^T\varphi_2^{-1}
\\
\implies \qquad s_2 & =  -u_2 -\frac{n}{2}\hat{\con}^{-1}_{\psi} - \frac{n}{2}\psi\psi^\T.
\end{align*}
\\
In much the same way that $s_{\psi}$ was found we find the maximum likelihood estimate for $\varphi$ as we vary $t$ in $s(t)$. The partial derivatives of $\ell\{\varphi, s(t)\}$ with respect to $\varphi$  are
\begin{align*}
\frac{\partial \ell}{\partial \varphi_1} & = \{u_1 + s_1(t)\} - n\varphi_2^{-1}\varphi_1,
\\
\implies \qquad \hat{\varphi}_2^{-1}(t)\hat{\varphi}_1(t) & = \frac{1}{n}\{u_1 + s_1(t)\} 
\\
& = \frac{1}{n}\{n\psi + t(u_1 - n\psi) \}
\\
& = \psi + t(\bar{y} - \psi).
\end{align*}
Notice that the above MLE agrees with what one might reasonably expect since $\varphi_2^{-1}\varphi_1$ is the mean vector. We use the above formula to solve for $\hat{\varphi}_2^{-1}(t)$: 
\begin{align*}
\frac{\partial \ell}{\partial \varphi_2} = & \{u_2 + s_2(t)\} + \frac{n}{2}\varphi_2^{-1} + \frac{n}{2}\varphi_2^{-1}\varphi_1\varphi_1^\T \varphi_2^{-1},
\\
\implies \qquad \hat{\varphi}_2^{-1}(t) = & \frac{2}{n} [-u_2 - s_2(t) - \frac{n}{2}\{\hat{\varphi}_2^{-1}(t)\hat{\varphi}_1(t)\}\{\hat{\varphi}_2^{-1}(t)\hat{\varphi}_1(t)\}^\T]
\\
 = & \frac{2}{n}\{ -u_2 + (1-t)(u_2 + \frac{n}{2}\hat{\con}_{\psi}^{-1} +  \frac{n}{2}\psi\psi^\T)\} -\frac{1}{n^2}\{u_1 + s_1(t)\}\{u_1 + s_1(t)\}^\T
\\
 = & (1-t)(\hat{\con}_{\psi}^{-1} + \psi\psi^\T) - \frac{2}{n}t u_2 -\{\psi + t(\bar{y} - \psi)\}\{\psi + t(\bar{y} - \psi)\}^\T
\\
 = & (1-t)\hat{\con}_{\psi}^{-1} + t(\psi\psi^\T - \bar{y}\psi^\T - \psi\bar{y}^\T)- \frac{2}{n}tu_2 - t^2(\bar{y} - \psi)(\bar{y} - \psi)^\T 
\\
= & \hat{\con}_{\psi}^{-1} - t^2(\bar{y} - \psi)(\bar{y} - \psi)^\T.
\end{align*}
As a check on our work we see that $\hat{\varphi}(0)$ provides the constrained MLE while $\hat{\varphi}(1)$ gives the unconstrained MLE. We see that $\hat{\varphi}_2^{-1}(t)$ is symmetric and thus $\hat{\varphi}_2(t)$ is symmetric meaning that any transpositions of these terms may be ignored in future calculations. Throughout we will use $\hat{\mu}(t)$ to represent $\hat{\varphi}_2^{-1}(t)\hat{\varphi}_1(t)$. Remembering that any terms not involving $t$ can be dropped from our calculations we find $\ell\{\hat{\varphi}(t);s(t)\}$ to be
\begin{align*}
\ell\{\hat{\varphi}(t);s(t)\}  = & \hat{\varphi}^\T(t) \{u + s(t)\} + \frac{n}{2}\log(\vert \hat{\varphi_2}(t) \vert) - \frac{n}{2}\hat{\varphi}_1^\T(t) \hat{\varphi}
_2^{-1}(t) \hat{\varphi}
_1(t)
\\
= & n\hat{\varphi}_1^\T(t) \hat{\mu}(t) + \text{Tr}\big[ \hat{\varphi}_2^\T(t) \{ u_2 + s_2(t)\} \big] - \frac{n}{2}\hat{\mu}^\T(t) \hat{\varphi}_2(t)\hat{\mu}(t) - \frac{n}{2}\log(\vert \hat{\varphi}_2^{-1}(t) \vert)
\\
= & \text{Tr}\big[\hat{\varphi}_2(t) \big\{ \frac{n}{2}\hat{\mu}(t)\hat{\mu}^\T(t) - \frac{n}{2}\con^{-1}_{\psi} - \frac{n}{2}\psi\psi^\T + t(u_2 + \frac{n}{2}\con^{-1}_{\psi}  + \frac{n}{2}\psi\psi^\T) \big\} \big] 
\\ & - \frac{n}{2}\log(\vert \hat{\varphi}_2^{-1}(t) \vert)     \\ = & \text{Tr}\big[ \hat{\varphi}_2(t)\big\{-\frac{n}{2}\con^{-1}_{\psi} + \frac{n}{2} t^2(\bar{y} - \psi)(\bar{y} - \psi)^\T\big\}\big] - \frac{n}{2}\log(\vert \hat{\varphi}_2^{-1}(t) \vert)
\\
= & \text{Tr}(-\frac{n}{2}I_p) - \frac{n}{2}\log(\vert \hat{\varphi}_2^{-1}(t) \vert)
\\
\equiv & - \frac{n}{2}\log(\vert \hat{\varphi}_2^{-1}(t) \vert).
\end{align*}
Similarly we find $\ell\{\hat{\varphi}(0);s(t)\}$ as
\begin{align*}
\ell\{\hat{\varphi}(0);s(t)\} = & \hat{\varphi}^\T(0) \{u + s(t)\} 
\\
= & n\hat{\mu}^\T(0) \hat{\varphi}_2(0)\hat{\mu}(t) + \text{Tr}[ \hat{\varphi}_2(0)\{u_2 + s_2(t)\}]
\\
= & \text{Tr}\big[ \hat{\varphi}_2(0)\big\{ n\hat{\mu}(t)\psi^\T + u_2 + s_2(t)\big\}\big]
\\
= & \text{Tr}\big[ \hat{\varphi}_2(0) \big\{ n\hat{\mu}(t)\psi^\T - \frac{n}{2}\con^{-1}_{\psi} - \frac{n}{2}\psi\psi^\T + t(u_2 + \frac{n}{2}\con^{-1}_{\psi}  + \frac{n}{2}\psi\psi^\T)\big\} \big]
\\
\equiv & \text{Tr}\big[ \hat{\varphi}_2(0) \big\{ tn(\bar{y} - \psi)\psi^\T + t(u_2 + \frac{n}{2}\con^{-1}_{\psi}  + \frac{n}{2}\psi\psi^\T)\big\} \big]
\\
= &  \text{Tr}\big[ \hat{\varphi}_2(0) \big\{ tn(\bar{y} - \psi)\psi^\T + t(n\psi\psi^\T - \frac{n}{2}\bar{y}\psi^\T - \frac{n}{2}\psi\bar{y}^\T)\big\} \big]
\\
= & 0.
\end{align*}
From here we can find the first piece of the conditional density:
\begin{align*}
\exp[\ell\{\hat{\varphi}(0);s(t)\}- \ell\{\hat{\varphi}(t);s(t)\}] \equiv \vert \hat{\varphi}_2^{-1}(t) \vert^{\frac{n}{2}}.
\end{align*}

\subsection*{S.6.2 Finding $\vert J_{\varphi\varphi}\{\hat{\varphi}(t) ; s(t)\}\vert$ and the nuisance parameter adjustment}

The likelihood in this scenario is the same as that in Example 5.3 of Fraser et. al. (2014). The canonical parameterization is also unchanged and so we can borrow the result that $$\vert J_{\varphi\varphi}\{\hat{\varphi}(t) ; s(t)\}\vert = \vert \hat{\varphi}_2^{-1}(t) \vert^{p + 2}.$$
The only term involving $t$ in the nuisance parameter adjustment is $\vert \ell_{\lambda\lambda}\{\hat{\varphi}(0);s(t)\} \vert$ which in turn only involves $t$ through $ \partial^2 / \partial \lambda^2 \{s^\T(t) \varphi\}$. Now $\con \mu$ is linear in $\lambda = \con$ and of course so is $\con$. As a result, all second order derivatives of $\varphi$ with respect to $\lambda$ disappear, meaning that the nuisance parameter adjustment is constant.

\subsection* {S.6.3 Making a change of variables}

Let's call $A = \hat{\con}_{\psi}^{-1}$, $v = (\bar{y} - \psi)$ and $B = \frac{1}{n}\sumonn (y_i - \bar{y})(y_i - \bar{y})^\T$. By definition $\vert \con^{-1}(t) \vert = \vert A - t^2vv^\T \vert$. Using the matrix determinant lemma we see that 
$$\vert \con^{-1}(t) \vert = \det(A) (1 - t^2 v^\T A^{-1} v) $$
By the constraint that $ \hat{\con}^{-1}(t)$ must be a valid covariance matrix, $t_{max}$ is the largest value of $t$ such that $\vert \hat{\con}^{-1}(t) \vert $ is positive definite. $\hat{\con}^{-1}(t)$ stops being positive definite as soon as one of its eigenvalues becomes $0$. Thus, $t_{max}$ is the solution to $\vert \hat{\con}^{-1}(t) \vert  = 0$:
\begin{align*}
\implies t_{max} = (v^\T A^{-1}v)^{-1/2}.
\end{align*}
Next we derive a formula that will be useful later on. It is easily seen that $A = B + vv^\T$. Using the Sherman-Morrison formula on $(B + vv^\T)^{-1}$ we see that 
\begin{align*}
v^\T A^{-1}v & = v^\T \bigg( B^{-1} - \frac{B^{-1}vv^\T B^{-1} }{1 + v^\T B^{-1}v} \bigg)v
\\
& = \frac{v^\T B^{-1}v}{1 + v^\T B^{-1}v}.
\end{align*}
Therefore, $(v^\T A^{-1}v) / (1 - v^\T A^{-1}v)   = v^\T B^{-1}v$. 
\\
\\
The directional $p$-value is given by 
\begin{align*}
p(\psi) = \frac{\int^{(v^\T A^{-1}v)^{-1/2}}_1t^{p-1} (1- t^2v^\T A^{-1}v)^{\frac{n - p - 2}{2}}dt}{\int^{(v^\T A^{-1}v)^{-1/2}}_0 t^{p-1} (1- t^2v^\T A^{-1}v)^{\frac{n - p - 2}{2}}dt}.
\end{align*}
For simplicity let $C = v^\T A^{-1}v$. The changes of variables made here are essentially identical to that in the normal linear regression example in Section 3. Make the change of variables $x = C t^2$ in the integral in the numerator:

\begin{align*}
\int^{C^{-1/2}}_1t^{p-1} (1- t^2C)^{\frac{n - p - 2}{2}}dt& = k\int^1_{C} x^{\frac{p-2}{2}}(1-x)^{\frac{n-p-2}{2}}dx.
\end{align*}

\noindent Make the change of variables to $x = 1/z$:

\begin{align*}
& = k\int^1_{C^{-1}}-\frac{1}{z^2} \big(\frac{1}{z}\big)^{\frac{p-2}{2}}\big(1-\frac{1}{z}\big)^{\frac{n-p-2}{2}}dz 
\\
& = k\int^{C^{-1}}_1 \big(\frac{1}{z}\big)^{\frac{p+2}{2}}\big(1-\frac{1}{z}\big)^{\frac{n-p-2}{2}}dz 
\\
& = k\int^{C^{-1}}_1 \big(\frac{1}{z}\big)^{\frac{n}{2}}(z-1)^{\frac{n-p-2}{2}}dz.
\end{align*}

\noindent Make the change of variables $t = z - 1$:

\begin{align*}
& = k\int^{C^{-1}-1}_0 \big(\frac{1}{t+1}\big)^{\frac{n}{2}}t^{\frac{n-p-2}{2}}dt 
\\
& = k\int^{\frac{1 - C}{C}}_0 \big(\frac{1}{t+1}\big)^{(\frac{n-p}{2}+\frac{p}{2})}t^{\frac{n-p}{2}-1}dt.
\end{align*}

\noindent Make the final change of variables $t = \frac{n-p}{p}x$:

\begin{align*}
& = k^\prime \int^{\frac{(1 - C)p}{C(n-p)}}_0 \bigg(\frac{1}{\frac{n-p}{p}x+1}\bigg)^{(\frac{n-p}{2}+\frac{p}{2})}\big(\frac{n-p}{p}x\big)^{\frac{n-p}{2}-1}dx
\\
& = k^{\prime\prime}\int^{\frac{(1 - C)p}{C(n-p)}}_0 \frac{\Gamma(\frac{n-p}{2} + \frac{p}{2})}{\Gamma(\frac{n-p}{2})\Gamma(\frac{p}{2})}\frac{n-p}{p-1} \big(\frac{1}{\frac{n-p}{p}x+1}\big)^{(\frac{n-p}{2}+\frac{p}{2})}\big(\frac{n-p}{p}x\big)^{\frac{n-p}{2}-1}dx
\\
& = k^{\prime\prime}P_W\bigg(W > \frac{n - p}{p}\frac{C}{1 - C}\bigg)
\\
& = k^{\prime\prime}P_W\bigg(W > \frac{n - p}{p}v^\T B^{-1}v\bigg),
\end{align*}

\noindent where $W \sim F(p,n - p)$. Now performing the exact same sequence of changes of variables on the integral in the denominator will result in a similar expression. The bounds of the integral will change as follows:
\begin{align*}
\int_{0}^{C^{-1/2}} \Rightarrow \int_0^1 \Rightarrow \int_1^\infty \Rightarrow \int_0^\infty.
\end{align*}

\noindent Thus the integral in the denominator of the directional $p$-value will equal $k^{\prime\prime}P_W(W > 0) = k^{\prime\prime}$. Hotelling's $T^2$ statistic is given by $T^2 = (n - 1)v^\T B^{-1}v$ and since $\frac{n - p}{p (n - 1)}T^2 \sim F_{p,n-p}$, the directional test is identical to the $p$-value obtained from Hotelling's $T^2$ test.

\end{document}